\documentclass[12pt]{amsart}

\usepackage{epsfig}

\newtheorem{theorem}{Theorem}[section]
\newtheorem{lemma}[theorem]{Lemma}
\newtheorem{cor}[theorem]{Corollary}
\newtheorem{proposition}[theorem]{Proposition}

\newtheorem{definition}[theorem]{Definition}

\newcommand{\Z}{\mathbb Z}

\newcommand{\al}{\alpha}

\newcommand{\comment}[1]{}
\newcommand{\spa}{{\rm span}}
\newcommand{\ov}{\overline}

\title{A Note on the Bar-Natan Skein Module}

\author{Marta Asaeda}

\address{Department of Mathematics, University of California Riverside,  900 Big Springs Drive, Riverside, CA, 92521 ,
 USA} 

\email{\tt marta@math.ucr.edu}

\author{Charles Frohman}

\address{Department of Mathematics, University of Iowa, Iowa City, IA
52242, USA}

\email{\tt frohman@math.uiowa.edu}

\thanks{The authors were sponsored in part by NSF grants
  \#DMS-0504199  and  \#DMS-0508635 respectively.}
\begin{document}
 
 \begin{abstract}
 We introduce a new skein module for three manifolds based on properly embedded surfaces and their relations introduced by D.Bar-Natan in \cite{BN1}, and modified by M.Khovanov \cite{M}. We compute the structure of the modules for some manifolds, including Seifert fibred manifolds. 
 \end{abstract}
 \maketitle
 
 \section{Introduction} This paper explores a new kind of skein module of three-manifolds that comes from Bar-Natan's study of the Khovanov homology of tangles and cobordisms \cite{BN1}. The module is a quotient of the free module on isotopy classes of {\em marked} surfaces in the underlying three-manifold. If the three-manifold is small in the sense that it has no incompressible surfaces, then the module is one dimensional on the empty surface. For sufficiently large manifolds its structure is more interesting. 
  
 The paper is organized as follows. In the next section we state the definitions and prove some propositions that will help to work examples. In section  3,  we look at two simple examples. Section 4 explores other normalizations for the module. Finally, we compute the Bar-Natan skein module for Seifert fibered spaces.
 
\section{The Basics}
\label{def}
Let $S$ be a surface. A simple closed curve $C$ embedded in $S$ is {\em essential}
if it does not bound a disk on $S$. A simple arc $A$ properly embedded in $S$
is {\em essential} if it does not cut a disk out of $S$. 

If $S$ is embedded in the
three-manifold $M$ a {\em compressing disk} for $S$ is a disk $D$ embedded in $M$ so
that $D\cap S=\partial D$ and $\partial D$ is essential on $S$. We say the
surface $S$ is {\em incompressible} if no component is a sphere bounding a ball and 
it has
no compressing disks.
% If $S$ is properly embedded in $M$ a boundary
%compression is a disk $D$ embedded in $M$ so that $D \cap S \subset \partial
%D$ is an essential arc $A$ in $S$ and $\partial D -A \subset \partial M$. We
%say that the surface $S\subset M$ is {\em boundary incompressible} if it does not
%have any components that are disks cutting a ball out of $M$ and $S$ admits no
%boundary compressions.  Finally, a properly embedded surface $S$ in $M$ is
%{\em essential} if it is incompressible and boundary incompressible.

A three-manifold $M$ is {\em irreducible} if every sphere bounds a ball. We
say the manifold $M$ is {\em prime} if every separating sphere in $M$ bounds a
ball.

\begin{definition} (Bar-Natan skein relations) \\
\begin{enumerate}
\item 
\raisebox{-.25in}{$\psfig{figure=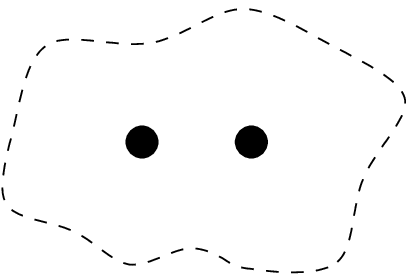,height=1.5cm}$} \ $=0$, 

\vspace{.25in}

\item  \raisebox{-.25in}{$ \psfig{figure=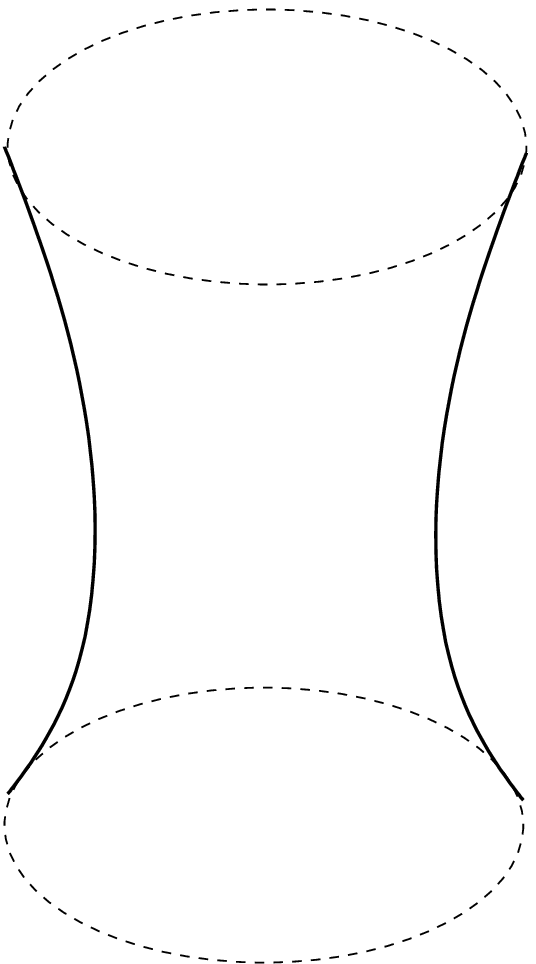,height=1.5cm}$} $
-$\raisebox{-.25in}{$ \psfig{figure=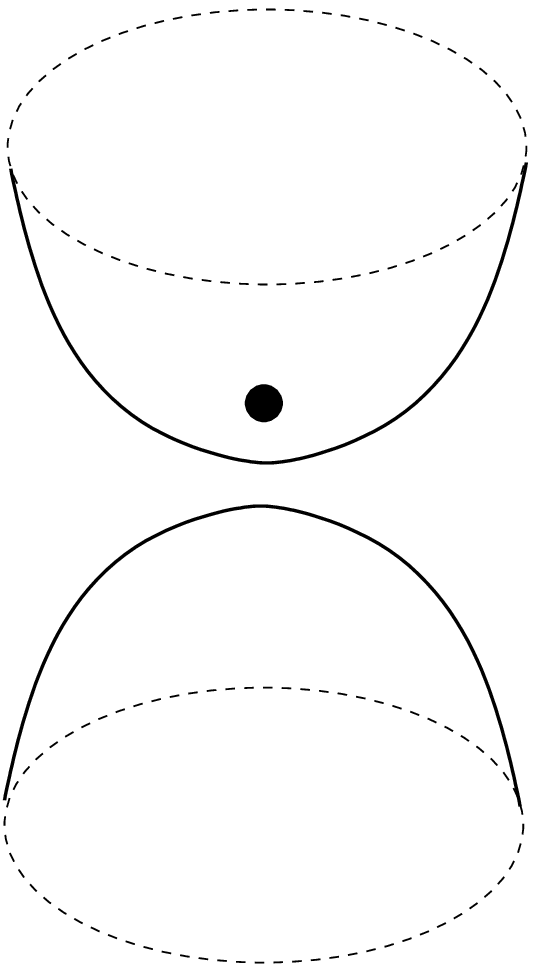,height=1.5cm}$} $
-$\raisebox{-.25in}{$\psfig{figure=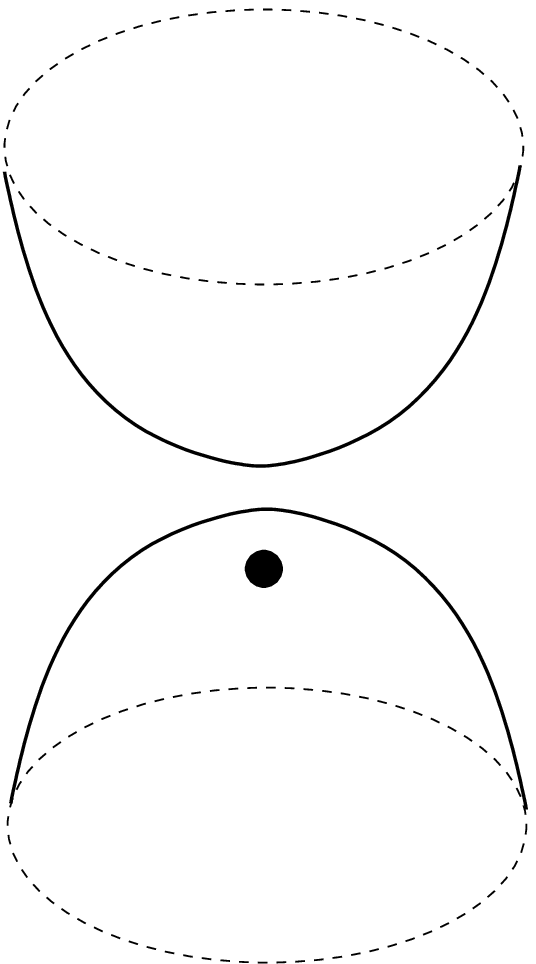,height=1.5cm} $}$=0,$

%\vspace{.25in} 
%\item  $2$ \raisebox{-.25in}{$ \psfig{figure=rel3left.eps,height=1.5cm}$}  $ -$
%\raisebox{-.25in}{$  \psfig{figure=rel3right.eps,height=1.5cm} $}$=0,$

\vspace{.25in}

\item \raisebox{-.13in}{$  \psfig{figure=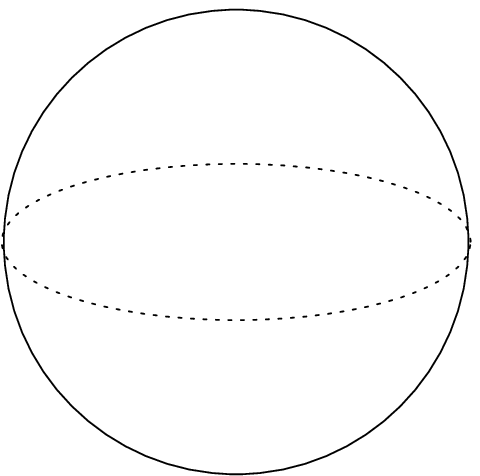,height=1cm}$}$ \sqcup F=0$ 

\vspace{.25in}

\item \raisebox{-.13in}{ $  \psfig{figure=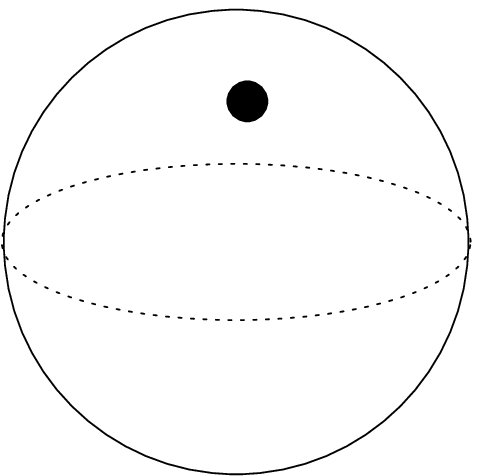,height=1cm}$}$ \sqcup F -F=0$ . 

\vspace{.25in}
In the last two diagrams we intend that the spheres bound balls. 

Any surface with more than one dot is zero under these relations. It is convenient to note 
 that the relation (2) implies \\
\vspace{.25in}
$2$ \raisebox{-.25in}{$ \psfig{figure=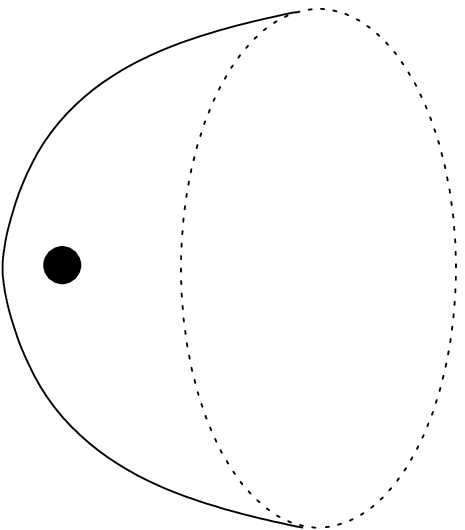,height=1.5cm}$}  $ -$
\raisebox{-.25in}{$  \psfig{figure=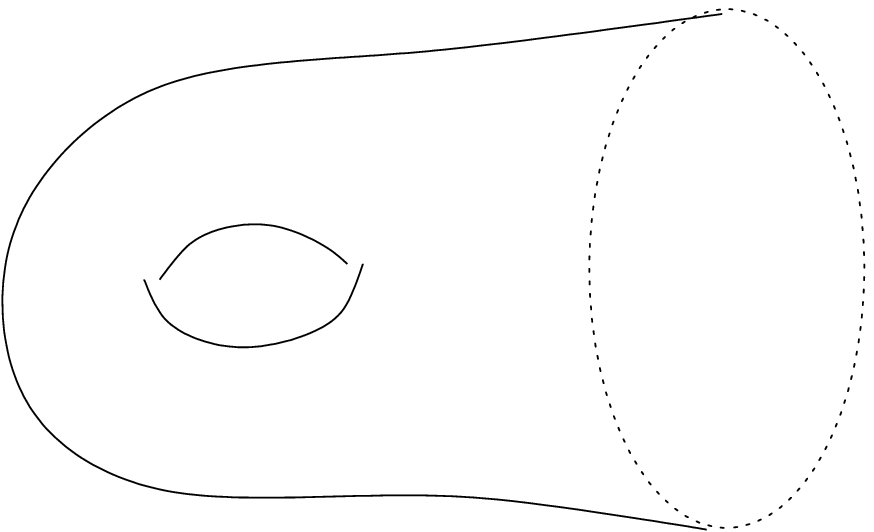,height=1.5cm} $}$=0.$
\vspace{.25in}

For 
later convenience, by "a dotted surface" we mean a surface with a dot on 
it, and "a white surface" we mean a surface without a dot.
\end{enumerate} 
\end{definition}
\begin{definition} (Bar-Natan skein modules) \\
%  \begin{itemize}
 % \item
    Let $M$ be a $3$-manifold  , and $R$ be a commutative unital ring.
A surface is marked if in addtion to the surface we indicate which components
carry  dots. Let $\mathcal{F}(M)$ be the set of isotopy classes of orientable marked 
surfaces, including the empty surface.    
Denote the free $R$-module with basis $\mathcal{F}$ by $R\mathcal{F}(M)$.
Let $S(M)$ denote the submodule of $R\mathcal{F}(M)$ spanned by the Bar-Natan
skein relations.  Let $BN(M,R)$ denote the quotient of $R\mathcal{F}(M)$ by
$S(M)$.   We call this the   Bar-Natan module of $M$

%
%\item Let $M$ be a $3$-manifold. Let $\mathcal{F}^o$ denote the set of isotopy
%  classes of oriented marked surfaces. If we perform the construction above,
%  except that we use  $R\mathcal{F}^o$ and skein relations between orientable 
%surfaces we denote the result $BN^o(M;R)$. This is the oriented Bar-Natan
%  skein module.

%\end{itemize}
\end{definition}
There are many variations on this definition depending on how we deal with orientations of surfaces. Note that one may define BN modules for a manifolds with boundary, using properly embedded surfaces. There are some variations for that as well, but we do not elaborate in that direction. 
%When we are making a statement that applies to both the oriented and
%unoriented Bar-Natan skein modules we will refer to it as $BN(M;R)$.

A {\em pure state} is an element of the Bar-Natan skein module that is
represented by a single marked surface. It is worth noting that a pure state
can have impure descriptions. We call a marked surface incompressible if the
underlying surface is incompressible.

\begin{proposition}
  $BN(M;R)$ is spanned by the pure incompressible states.
\end{proposition}

\proof Starting with a representation of an element  $BN(M,R)$  as a linear
combination of pure states. If one of the diagrams is compressible, then apply
the second skein relation to replace it by the sum of two diagrams after doing
the compression. We also eliminate spheres bounding balls using the fourth and fifth relations. After finitely many steps, the skein will be written as a sum
of incompressible pure states. \qed

It is evident  that $BN(S^3)=R\phi$, that is, it is the free module of
rank $1$ over $R$ with basis the empty surface.  How does a given surface in
the three-sphere evaluate?  The skein relations tell us how to deal with a
sphere. Suppose that the surface is a torus. A torus embedded in the
three-sphere always has a compressing disk. Suppose we have a simple state
consisting of the torus with no dots. After one compression, we have written
as a sum of two states consisting of a sphere with one dot, so the torus
evaluates as $2\emptyset$. If the torus has a dot, then after an application of the
second skein relation we have a sum of two spheres with two dots, which
evaluates to zero.  Using these rules it is easy to evaluate any surface in
the three-sphere.

We say two disjoint $S_0$ and $S_1$ in $M$ are parallel, 
if there is an embedding of $S\times [0,1]$ into $M$ so $S\times\{0\}$ is
$S_0$ and $S_1$ is $S \times \{1\}$. We call image of the embedding of
$S\times [0,1]$ the region between them.

\begin{proposition}
\label{parallel}
  Suppose a surface $F$ contains two parallel connected surfaces $S_0,  S_1$  so
  that no other surface in the state touches the region between them. 
Let $P_0$ be a pure state based on $F$  so that $S_0$ carries a dot and $S_1$
  does not carry a dot. Let $P_1$ be the pure state that is the same away from
  $S_0$ and $S_1$, but now $S_1$ carries a dot and $S_0$ does not carry a dot.
Then $P_0=-P_1$. 

Suppose that  $T$ is a pure state so that both $S_0$ and $S_1$ carry a
dot then $T$.  If $S_0$ is a sphere, then $T$ is equivalent to the pure state obtained by deleting $S_0 \cup S_1$.  If $S_0$ has genus greater than zero then $T=0$.
\end{proposition}

\proof   First we assume that $F=S_0\cup S_1$. By $A^{\bullet}$ we denote a connected surface $A$ carrying a dot. Then we may have an expression of the pure state $P_0= {S_0}^{\bullet} \cup S_1$, $P_1= {S_0}\cup {S_1}^{\bullet}$ respectively. By the second skein relation we have 
${S_0}^{\bullet} \cup S_1+{S_0}\cup {S_1}^{\bullet}=S_0 \sharp S_1 $, where the connected sum is taken through the region between $S_i$'s. Note that the right hand side bounds a handle body $S \backslash \{disk\} \times [0,1]$. which has the genus $0$ or greater than $1$.  Therefore it determines zero state, thus we obtain $P_0 = -P_1$. Similarly $T={S_0}^{\bullet}  \cup {S_1}^{\bullet} =(S_0 \sharp S_1)^\bullet - {S_0}  \cup {S_1}^{\bullet \bullet} =(S_0 \sharp S_1)^\bullet. $ If $S$ is sphere $S_0 \sharp S_1$ is a sphere bounding a ball, thus we obtain $T=1$, otherwise $T=0$. Since all the procedures are done within the region between $S_i$'s, the general case where $F$ has more components than $S_i$'s is proved in the same way. \qed. \\

In general you cannot move dots around so fluidly. To this end we can define a
filtration of $BN(M,R)$. We say a surface  has level $m$ if it has $m$ connected
components.  Let $F_m$ be the span in $BN(M,R)$ of all surfaces having level
less than or equal to $m$. It is clear that $\cup_m F_m=BN(M,R)$. Let
$G_m=F_m/F_{m-1}$  where $F_{-1}=\{0\}$.  By the graded object we mean
$\sum_m G_m$. An important property of the graded object is:

\begin{proposition} \label{pure}
  If $P_0$ and $P_1$ are pure states as above of level $m$, then $P_0=-P_1$
in $G_m$. If a pure state of level $m$ has more than one dot, and
$M$ is connected then that state is zero in $G_m$.
\end{proposition}

\proof If you tube $S_0$ to $S_1$ you get a state of level $m-1$ so it is zero
in the quotient. That proves the first statement. To prove the second
statement use the first relation to shift dots till they are on the same
component. \qed

There is one more final proposition we would like to make here about the relations. 

\begin{proposition}
  \label{independence}
Suppose that $M$ is prime.
If
$P_i$ is a collection of distinct pure incompressible states without dots then
they are linearly independent. 
\end{proposition} 

\proof To prove this result we construct a family of linear functionals on
$\mathcal{F}(M)$ that is dual to the elements induced by the $P_i$ and
vanishes
on $S(M)$.  Let $F_i:BN(M,R)
\rightarrow R$ be defined as follows. Given a marked surface $S$, if
it is a surface isotopic to $P_i$ along with a disjoint union of
spheres bounding balls carrying a dot, and compressible tori that don't have a
dot, so that compressing any torus yields a sphere bounding a ball,  then let $F_i(S)=(1/2)^k$ where $k$ is the number of compressible tori without a
dot. Otherwise, let $F_i(S)=0$.

If this map were not well defined it would be because there are two ways to
compress a particular torus. One that yields a sphere bounding a ball, and one
that yields a sphere that does not bound a ball. If this happens the two
compressing disks must lie on opposite sides of the torus, which is
separating. The side containing the disk, so that compressing along the disk
yields a ball, is a solid torus. The other side can be seen as the connect sum
of a solid torus with a $3$-manifold that is not a sphere. Hence if the two
disks like this exist, $M$ is the connected sum of a nontrivial manifold with
a lens space, and is not prime. Therefore $F_i$ is well defined. 
It is easy to check that $F$ sends all the relations in $BN(M,R)$ to zero, so it
induces a linear functional on $BN(M)$. The family of linear functionals $F_i$
is dual to the $P_i$ in the sense that $F_i(P_j)=\delta_i^j$. Therefore the
$P_i$ are independent.

\section{Two Examples}

From the results of the last section it is evident that if $M$ has no
incompressible surfaces, $BN(M,R)$ is free of rank one on the empty
surface. To get an interesting answer, you need to have plenty of
incompressible surfaces in your three-manifold. In this section we will discuss
two examples. The first is $S^1\times S^2$. The only incompressible surfaces
in $S^1\times S^2$ are surfaces whose connected components are all isotopic to
$\{*\}\times S^2$. Hence the collection is nonempty but has a really simple
structure. The second example we discuss is $S^1\times S^1\times
S^1$. 
%Although the collection of incompressible surfaces is much more
%complicated 
The answer turns out to be very similar to the first example.

Let $M$ be $S^1 \times S^2$. We compute $BN(M,R)$. Assume that $2$ is
invertible in $R$.
\begin{theorem} 
By $z^k$ we denote a disjoint union of $k$ parallel spheres in $M$ that are
homologically non-trivial. 
Let $e_m$ be a disjoint union of $z^m$ and a dotted non-trivial sphere
(i.e. $m+1$ parallel 
non-trivial spheres with a dot on one of them). Then
$$BN(M,R)= R[z]  \oplus R e_0. $$

\end{theorem}

\proof There are two statements that need to be proved. The surfaces span, and
the surfaces are independent. in the following we only consider pure states based on parallel copies of essential $S^2$ unless otherwise is specified. \\

 The lemmata below follow from Proposition \ref{parallel}.  
\begin{lemma} (Shifting)
\label{switch}
Let $\alpha$ be a dotted sphere  in $e$ and $\beta$ be a white sphere parallel to $\alpha$. Then $e=-e''$, where $e''$ is obtained from $e$ by moving the dot from $\alpha$ to $\beta$. 
\end{lemma}
 \begin{lemma} (reduction)
\label{reduce}
If $e$ contains two parallel dotted spheres $\alpha_1$, $\alpha_2$, then $e=e'$, where $ e'=e-  \{\alpha_1, \alpha_2\}$. 
\end{lemma}
 Using these lemmata we may shift dots and reduce parallel dotted spheres. Thus one may describe a pure state  with more than one dots as a linear combination of the states with no dot, or a dot on some component.  However, a pure state with a dot is zero if the number of component is greater than one: while shifting of a dot from a sphere to an adjacent sphere causes sign change,   it is on the other hand isotopic to the original state, thus the same element in the module. 
%  {\bf end of insert:jump to Charlie's independence argument starting "Next we need to see". } 
%  
%Suppose that a surface consists of some parallel
%copies of nontrivial spheres so that one of them carries a dot, and there are
%at least two connected components.  
%If two adjacent surfaces carry a dot, we can tube two adjacent spheres together
%with one dot. By the second skein relation this  new surface is equivalent to
%the first surface. The tubed together spheres bound a ball and it has a dot,
%so we may eliminate the component and the resulting skein has two fewer
%components.

%If two adjacent
%components have the property that one carries a dot and the other doesn't,
%then using the tubing trick above we can move the dot from one to the other
%and the affect is the to send the skein to its negative. If the new skein has
%two adjacent components carrying a dot then apply the first construction.  If
%the surface has more than one dot, we can keep permuting dots till they are
%adjacent and the apply the first construction. If not, then after one
%permutaion of dots, do an isotopy to take the surface back to itself. Hence
%twice the surface is equal to $0$ in $BN(M)$. However, we assume that $2$ is
%invertible in $R$ which implies the skein induced by the original surface is
%zero.

%Hence if a surface has dots we can write it as a linear combination of
%incompressible surfaces without dots, and a single incompressible sphere with
%a dot. 

Next we need to see that our skeins are independent.
 By Proposition
\ref{pure} the skeins $z^k$ are all independent. We just need to show that
$e_0$ is independent from the $z^k$. To this end we define a linear functional
$E:R\mathcal{F}({S^1\times S^2}) \rightarrow R$.  If a surface is  the result of
taking the disjoint union a family of spheres with dots  bounding balls and
compressible tori without a dot that are nullhomologous in $H_2(S^1\times
S^2;\Z_2)$, with $e_0$ or a torus without a dot that represents something
nonzero in $H_2(S^1\times
S^2;\Z_2)$  send it  to $(1/2)^k$ where $k$ is the
number of torus components.  Send any other surface to $0$.  Once again it is easy to check that this linear
functional vanishes on all the relations so descends to $BN(S^1\times
S^2,R)$. Furthermore, it sends $e_0$ to $1$ and all $z^k$ to zero. 
Ergo the surface $e_0$ is independent from the surfaces $z^k$.\qed

The case of $T^3=S^1\times S^1\times S^1$ is so similar we will just posit the
answer without proof.  The nonempty connected orientable surfaces in $T^3$ are all
tori. Up to istopy they are in on to one correspondence with the pairs of triples of integers $\pm(p,q,r)$ that
are relatively prime and not all zero. The module $BN(T^3)$ has a basis consisting of the empty surface,
parallel copies of the tori paramatrized above without a dot, and single
copies of the tori above carrying one dot.

\section{An exploration of normalization}
 In this section we investigate if there is any other choice for the constants in Bar-Natan skein relations as in Definition \ref{def}.  
 First let 
 $$  \raisebox{-.13in}{\psfig{figure=rel4.eps,height=1cm}} =y,$$
 $$ \raisebox{-.13in}{ \psfig{figure=rel5.eps,height=1cm}} =z,$$
 $$ \raisebox{-.13in}{\psfig{figure=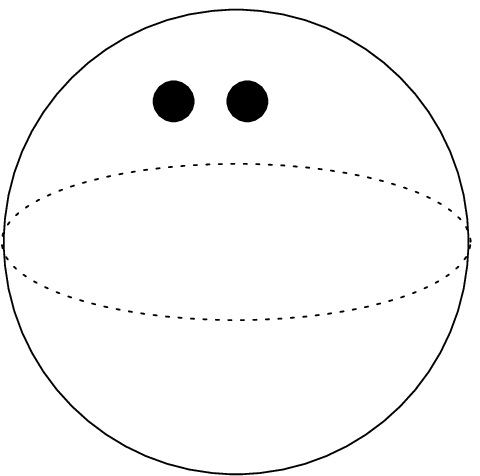,height=1cm}} =x. $$
By using the relation (2) in Definition \ref{def} we have
 $$  \raisebox{-.13in}{\psfig{figure=rel4.eps,height=1cm}}= y =
\raisebox{-.13in}{\psfig{figure=rel4.eps,height=1cm}\
  \psfig{figure=rel5.eps,height=1cm}}+
\raisebox{-.13in}{\psfig{figure=rel5.eps,height=1cm}\ \psfig{figure=rel4.eps,height=1cm}}=2yz$$
 and
 $$\raisebox{-.13in}{\psfig{figure=rel5.eps,height=1cm}}= 
\raisebox{-.13in}{\psfig{figure=rel5.eps,height=1cm}\
  \psfig{figure=rel5.eps,height=1cm}} +
\raisebox{-.13in}{\psfig{figure=rel4.eps,height=1cm}\  \psfig{figure=rel6mod.eps,height=1cm}}.$$
 we need $2,y$ to be invertible in $R$ if $y\neq 0$, then we have $z=1/2$, and $xy=1/4$. In the case $y=0$, we have $z=z^2$, then $z=0$ or $1$. For $z=0$ we have totally trivial theory where everything is zero, and the latter case we have the relation identical to the Bar-Natan skein relation. Thus we have the following new possible  skein relations:
 \begin{definition} (modified Bar-Natan skein relations) \\
%relation 1: 
%$\psfig{figure=rel1.eps,height=1.5cm} =0$, 
%\\
\begin{enumerate}

\item  \raisebox{-.25in}{$ \psfig{figure=rel2left.eps,height=1.5cm}$} $
-$\raisebox{-.25in}{$ \psfig{figure=rel2right1.eps,height=1.5cm}$} $
-$\raisebox{-.25in}{$\psfig{figure=rel2right2.eps,height=1.5cm} $}$=0,$
 %
%\item relation 3: $2\  \raisebox{-.25in}{\psfig{figure=rel3left.eps,height=1.5cm}}
%=   \raisebox{-.25in}{\psfig{figure=rel3right.eps,height=1.5cm}} $  

\vspace{.25in}

\item  $   \raisebox{-.13in}{\psfig{figure=rel4.eps,height=1cm}} -y=0$ 
\vspace{.25in}

\item $   \raisebox{-.13in}{\psfig{figure=rel5.eps,height=1cm}} -1/2=0$  
\vspace{.25in}

\item  $  \raisebox{-.13in}{\psfig{figure=rel6mod.eps,height=1cm}} -x=0,$ \\
\vspace{.25in}
where $y=\frac{1}{4x}$.
 \end{enumerate}\end{definition}
 %%%%%%%%
 These relations turn out to imply that the distribution of the dots in fact do not matter as demonstrated next. Let $\al_k$ be a pure state with $k$ dots. Then \\
 \begin{eqnarray*}
\al_k= \raisebox{-.13in}{\psfig{figure=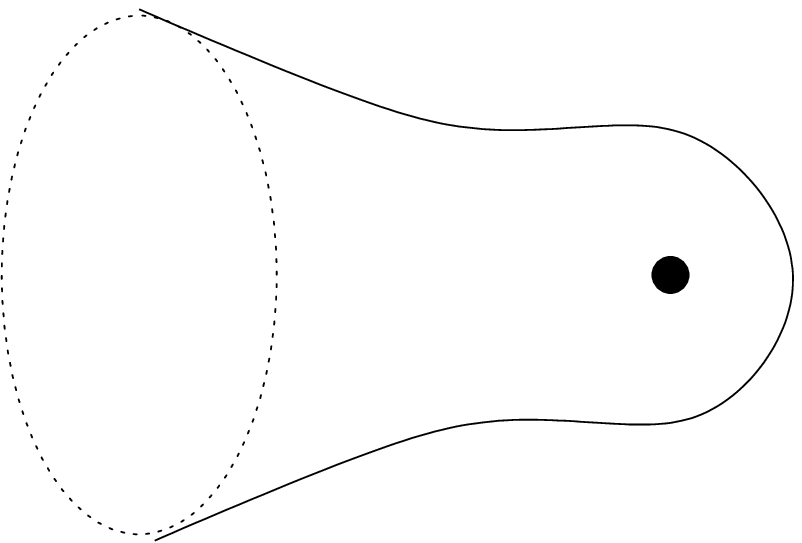,height=1cm}}&=&\raisebox{-.13in}{\psfig{figure=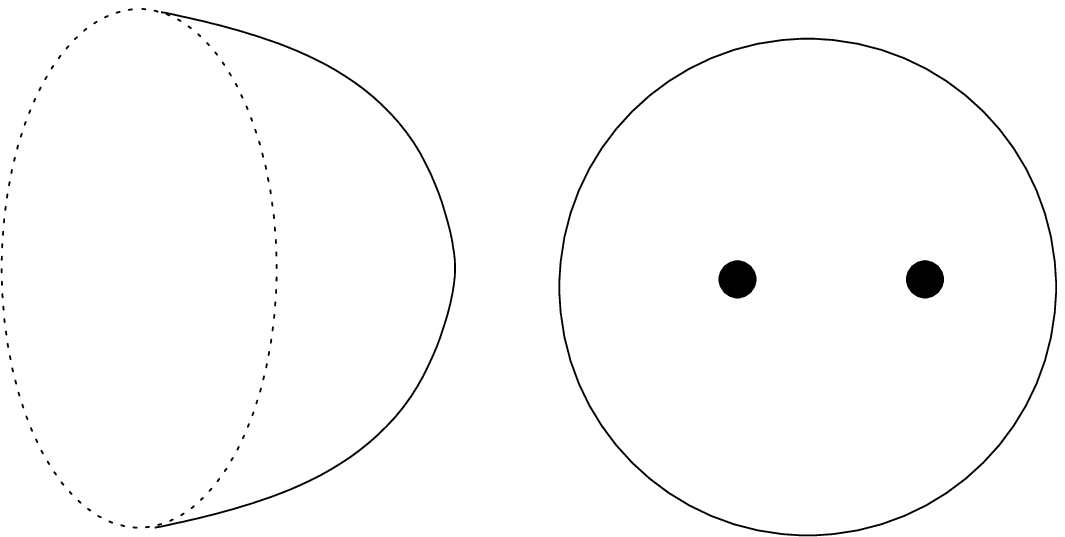,height=1cm}}+\raisebox{-.13in}{\psfig{figure=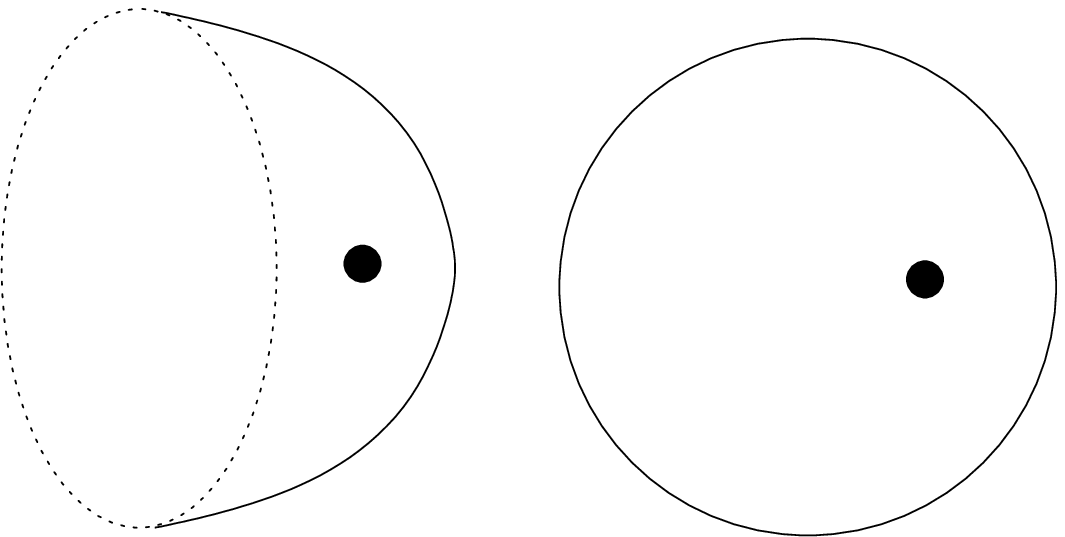,height=1cm}} \\
 &=& x \al_{k-1} + \frac{1}{2} \al_k,
 \end{eqnarray*}
 thus we have $\al_k=2x  \al_{k-1}= ... (2x)^k S_0$. Thus the first relation is simplified as \\
 \begin{enumerate}
\item[(1)'] \raisebox{-.25in}{$ \psfig{figure=rel2left.eps,height=1.5cm}$} $
-2$\raisebox{-.25in}{$ \psfig{figure=rel2right1.eps,height=1.5cm}$} $
=0,$
\end{enumerate}
   By scaling $\al_k \to 1/2^k \al_k=:\al'_k$, we get the following equivalent definition.  
    \begin{definition} (doubly modified Bar-Natan (MBN) skein relations) \\
%relation 1: 
%$\psfig{figure=rel1.eps,height=1.5cm} =0$, 
%\\
   \begin{enumerate}
   \item  $\raisebox{-.25in}{
   \psfig{figure=rel2left.eps,height=1.5cm}}\  =
   \raisebox{-.25in}{\psfig{figure=rel2right1.eps,height=1.5cm}}   $ \\ (the dot may go with the upper component or the lower one.) 
\vspace{.25in}
%relation 3: $ \psfig{figure=rel3left.eps,height=1.5cm}   =  \psfig{figure=rel3right.eps,height=1.5cm} $  \\
\item  $  \raisebox{-.13in}{\psfig{figure=rel4.eps,height=1cm}}
  =x^{-1}$ 

\item  $  \raisebox{-.13in}{\psfig{figure=rel5.eps,height=1cm}} =1$ . 

\item $\raisebox{-.13in}{ \psfig{figure=rel6mod.eps,height=1cm}} =x. $ 
\end{enumerate}
 \end{definition}
 Let $MBN(M,R):= $ \\ $ {\rm span}_R\{\mbox{\rm closed dotted or white surfaces in } M \} /
\mbox{\rm MBN relations}.$ \\
Note that in this definition we do not limit the surfaces to be orientable.  It is easy to check that $\al'_k=x\al'_{k-1}$ for a pure state $\al'_k$ with $k$ dots in $MBN(M,R)$, i.e. a dot corresponds to multiplication by $x$. Thus a compression also gives rise to multiplication by $x$. 
% Because of the compressing relation, $MBN(M,R)$ is generated by incompressible surfaces. 
 Taking the contribution of dots and white spheres bounding balls into account, we have 
    $$  MBN(M,R)  \cong R[x^\pm]H_2 (M, \Z_2) \cong R[x^\pm] \otimes H_2 (M, \Z_2)$$
with the identification of the elements given by $F \in MBN(M,R) \to x^{g(F)-|F|} F \in  R[x^\pm]H_2 (M, \Z_2)$, where $F$ is a surface without dots, $g(F)$ is the genus of $F$, and $|F|$ is the number of the connected components of $F$. 
Note that $MBN(M,R)$ has an algebra structure, equipped with the product induced by the double curve sum on the homology group. 

\section{Seifert fibred spaces}

 A Seifert fibered space is a compact three-manifold $M$ that is the disjoint union of
circles,
so that every circle has a neighborhood that is a union of circles, and 
the decomposition of the neighborhood is modeled on a standard 
$(p,q)$ decomposition of a solid
torus into circles (see \cite{He}), or if the circle lies in the boundary
of $M$, the neighborhood looks like half a solid torus fibered by longitudes.
 The decomposition space obtained by
identifying two points if they lie on the same circle is a surface $F$. The
quotient
map $\psi:M \rightarrow F$ has the property that its restriction to the
complement of finitely many fibers is an $S^1$ bundle. The fibers that must be
deleted are called the {\em singular} fibers. These correspond to points where
the local model is built on $(p,q)$ where $|p|>1$.  A surface $S$ in a Seifert
fibered space is called {\em vertical} if it is a union of fibers. It is
called {\em horizontal} if it is transverse to all the fibers. If $S$ is
horizontal in $M$ then the result of cutting $M$ along $S$ is homeomorphic to
$S \times [0,1]$. We will restrict our attention to Seifert fibered spaces
where the total space is orientable as is the decomposition space. In this
case it is a theorem of Waldhausen \cite{Wa} that every two-sided incompressible surface
is isotopic to a vertical surface or a horizontal surface.

 In the following we study the structure of $BN(M,R)$ where $M$ is a Seifert fibered space, $F$ is its orbifold, with $l$ singularities $p_1, ..., p_l$, and $\psi:M \to F$ be the fibration map always. We assume  that $F$ is closed, connected, and orientable for simplicity. By Waldhausen's theorem we immediately obtain the following decomposition:
 $$BN(M,R)=BN_v(M,R) \oplus BN_h(M,R),$$
 where $BN_v(M,R)$ is generated by vertical surfaces with dots, and $BN_h(M,R)$ is generated by horizontal surfaces with dots. 
%As discussed above, most of the orientable incompressible surfaces in $M$ are vertical. 

For now we focus on $BN_v(M,R)$. 
A vertical surface is   a pre-image of simple closed curves on $F$.  Connected sums of simple closed curves on $F$ correspond to "annular sums" of vertical surfaces along  longitudes: 
\begin{definition} (annular sum) 
Let $N$ be a three manifold. Let $S_1$, $S_2$ be connected components of a  surface properly embedded in $N$, and $A$ be an annulus embedded in $N$ so that each of two boundary component $c_i$ lies in $S_i$ respectively, and that $c_i$ is orientation preserving on $S_i$. We define the annular sum of $S_1$ and $S_2$ along with $A$ by 
$$ S_1 \sharp_A S_2 = S_1 \backslash N(c_1) \cup S_2 \backslash N(c_2) \cup A' \cup A'',$$ 
where $A'$ and $A''$ are parallel copies of $A$ so that for their boundary components $\partial A'=c'_1 \cup c'_2$ and $\partial A''=c''_1 \cup c''_2$, $c'_i \cup c''_i = \partial {\ov{N(c_i)}}$ holds. 
%
%Let $N$ be a three manifold. Let $S_1$, $S_2$ be connected components of a  surface properly embedded in $N$, containing orientation preserving simple closed curves $c_1$ and $c_2$. Suppose that there is an annulus $A$ embedded in $N$ so that $\partial A=A \cap (S_1 \cup S_2) = c_1 \cup c_2$. We define an annular sum
%$S_1 \sharp_A S_2$ by removing a regular neighborhood $A_i$ of $c_i$ from $S_i$ and identifying $\partial A_1$ and $\partial A_2$.   
\end{definition}
The above definition is illustrated as follows:  
\vspace{5mm}
$$S_1 \cup S_2 =\raisebox{-.25in}{\psfig{figure=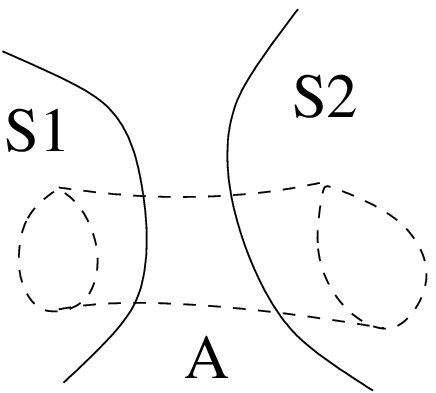,height=2cm}} \Rightarrow 
 S_1 \sharp_A S_2 =\raisebox{-.25in}{\psfig{figure=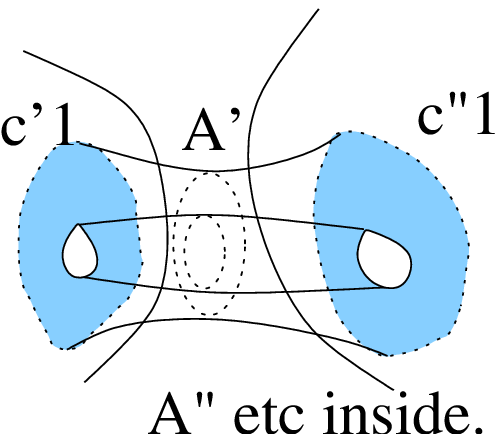,height=2cm}}$$
 \vspace{5mm}
 
Note that the annular sum depends on choice of  $A$.  One may also consider the case $S_1=S_2$ in the same way.  

We have the following skein relation associated with annular sum, derived from B-N skein relation. 
\begin{proposition}
\label{annskein}
For the above setting, let $S'_1$, $S'_2$ be the components of $S_1 \sharp_A S_2$ that contain components of $\partial A_i$ each. Then we have 
$$S_1^{\bullet} \cup S_2+ S_1 \cup S_2^{\bullet}  = S_1 \sharp S_2 = S'^{\bullet}_1 \cup S'_2+ S'_1 \cup S'^{\bullet}_2, $$
where by $S^\bullet$ we mean a singly dotted component $S$. 
\end{proposition}
This is easy to observe by the fact that both $S_1 \cup S_2$ and $S'_1 \cup S'_2$ are obtained by  applying the relation 2 in BN relations to  $S_1 \sharp S_2$ in different directions. 
Note that it is possible that $S_1=S_2$, or $S'_1=S'_2$, in those cases we have just  twice of the same term in each formula. Indeed by our condition that the boundary components of $A$ need to be orientation preserving,  if $S_1$ and $S_2$ are distinct components, $S'_i$ are in fact one component, and vice versa. 

Now we come back to our Seifert fibred manifold $M$. In the following we
consider the case where the surfaces are all vertical surfaces of $M$, and the annuli used for annular sum are also vertical. In the following we define a skein module on $F$ based on the curves corresponding to vertical surfaces, with relation induced by BN relations along with the relation given by Proposition \ref{annskein}. 
\begin{lemma}
A vertical surface is separating (resp. non-separating) in $M$ if and only if the corresponding simple closed curve on $F$ is separating (resp. non-separating). An annular sum of two vertical surfaces along the longitude corresponds to a band sum of corresponding s.s.c's in $F$. 
\end{lemma}
\begin{definition}
\label{surfaceskein}
For a surface $F$ we define
\begin{eqnarray*}
&& SBN(F,R) \\
&:=&{\rm span}_R\{ \mbox{\rm   properly embedded curves with or without dots}\} \slash relations, \end{eqnarray*}
where the relations are given by the following: \\ \ \\
{\bf SBN relations} \\ 
\begin{enumerate}
\item $\raisebox{-.25in}{\psfig{figure=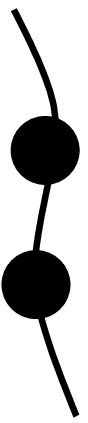,height=1.5cm}} =0$, \vspace{0.25in}
%relation 2:  $\psfig{figure=2curves1.eps,height=1cm} + \psfig{figure=2curves2.eps,height=1cm}=
 %2 \psfig{figure=glued.eps,height=1cm}, $ \\
\item  $\raisebox{-.15in}{\psfig{figure=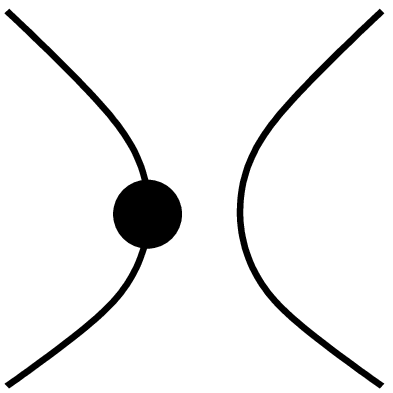,height=1cm}} + \raisebox{-.15in}{\psfig{figure=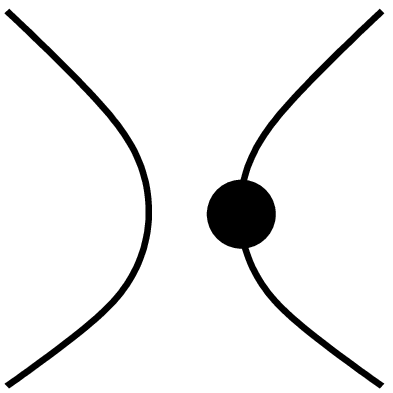,height=1cm}}=
 \raisebox{-.15in}{\psfig{figure=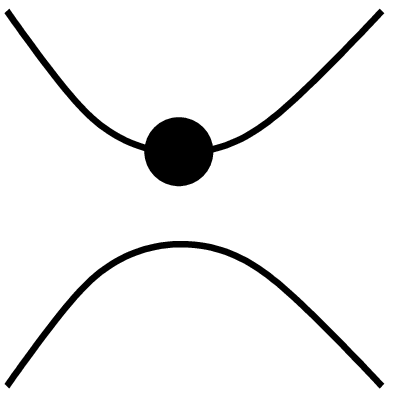,height=1cm}} + \raisebox{-.15in}{\psfig{figure=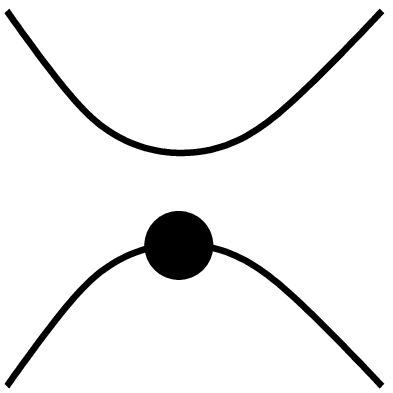,height=1cm}}$  \vspace{0.25in}
\item $\raisebox{-.15in}{\psfig{figure=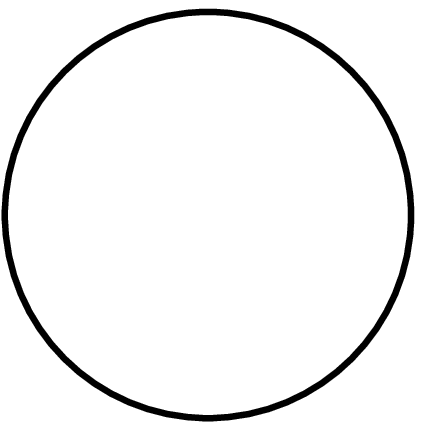,height=1cm}}=2, $  \vspace{0.25in}
\item $ \raisebox{-.15in}{\psfig{figure=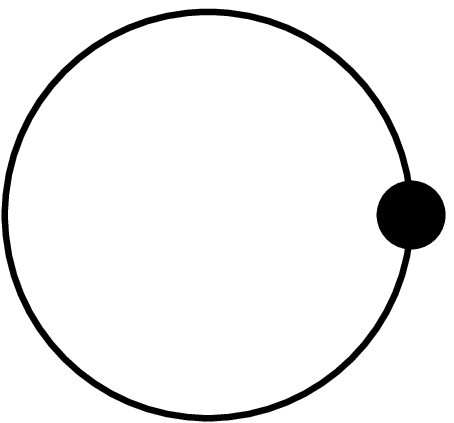,height=1cm} }=0,$  \vspace{0.25in}
\end{enumerate}
where  the circles   bound a disk, possibly with one singular point. 
\end{definition}
% \oplus \oplus_{z} \{z^k\}_k \oplus \{d_z\},$$
% where by $z^k$ we denote the $k$-parallel copies of a horizontal  surface $z$, and by $d_z$ we denote a dotted copy parallel to $z$. 
%We shall prove this proposition as well as the parametrization of horizontal surfaces after computing $SBN(F,R)$ in the following. 
By construction it is clear that $BN_v(M,R) \hookrightarrow SBN(F,R)$. In the following we study the structure of $SBN(F,R)$. 
\subsection{Structure of $SBN(F,R)$}  \ \\ \ \\
The SBN relation 2 is typically realized as the following situation for closed curves: 
\[ \raisebox{-.15in}{\psfig{figure=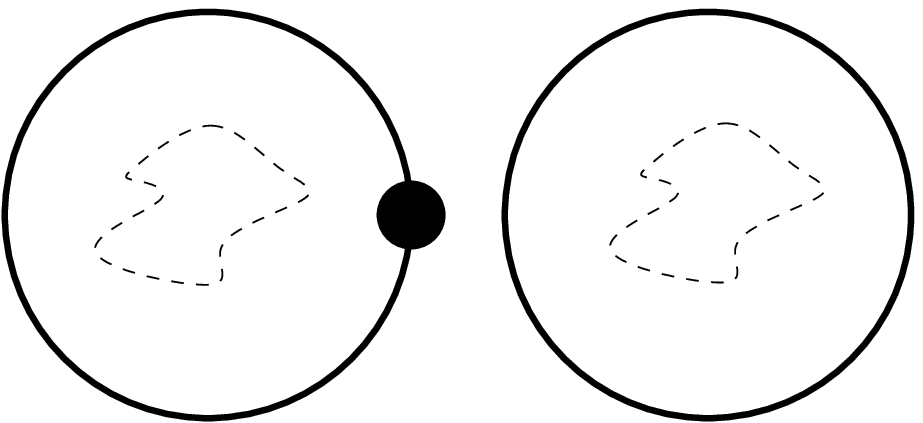,height=1cm}} +\raisebox{-.15in}{ \psfig{figure=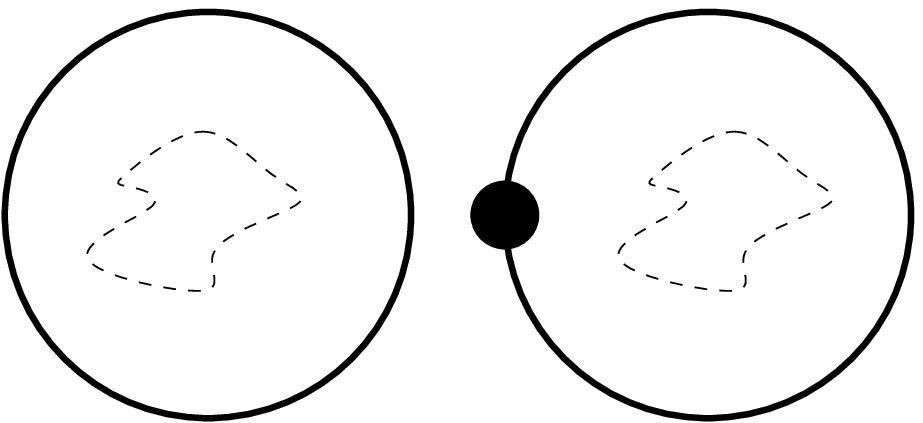,height=1cm}}=
 2\  \raisebox{-.15in}{\psfig{figure=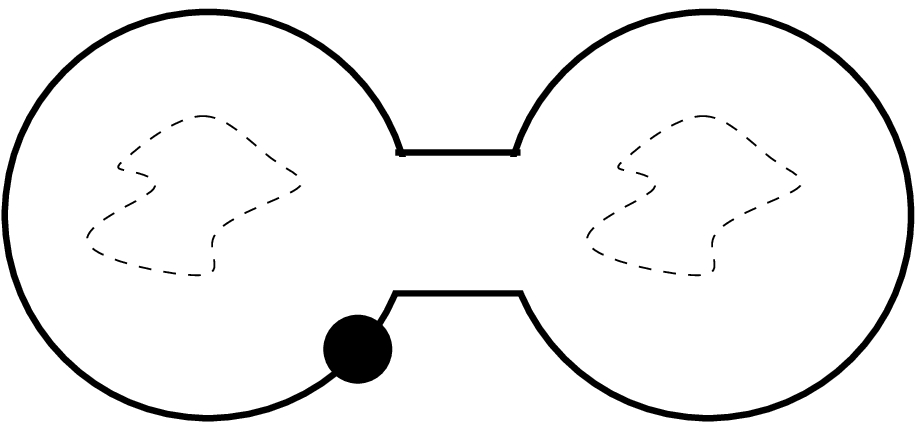,height=1cm}}.\] 
 %\\ where at least one of the two curves in the left hand side is orientation preserving on $F$, \\ \ \\
% (2) $2 \psfig{figure=orirevleft.eps,height=1.2cm} = 2 \psfig{figure=orirevright.eps,height=1.2cm},  $ \\where the curves are orientation reversing on $F$ (there is no crossing). \\ \ \\
%%%%%%%%%%%%%%%%%%%%
%Thus the understanding of $BN(M,R)$ is reduced to the understanding of $SBN(F,R)$, which may be possibly simpler. 
We call an element $\al \in SBN(F,R)$ pure state if it determines a pure state in $BN(M,R)$. All notions, such as levels, introduced for states in $BN(M,R)$ are used for corresponding states in $SBN(F,R)$ without explanation. 

Our goal is to find a basis for $SBN(F,R)$. Since we know from the earlier section that all the white pure states are linearly independent, we focus on describing the relations between dotted pure states. Since our focus is mainly on pure states, we will simply use the word "state" to mean a pure state unless otherwise is stated.  We provide some tools in the following. 
%\begin{proposition}
%BN(M,R) is generated by disjoint unions of pre-images of  rnsc's with or without dots, white vertical incompressible surfaces, and white or dotted horizontal surfaces. 
%\end{proposition}
%\begin{definition}
%We call $\al \in SBN(F,R)$ a pure state if it is expressed as a disjoint union of non-trivial circles. We call a circle or a curve  on $F$ white if it has no dot, and dotted, or black if it has a dot. 
%\end{definition}
\begin{proposition} 
 $2 c^\bullet=0$ if $c$ is a separating curve in $F$. 
 \end{proposition}
 This is from the following lemmas. 
\begin{lemma}
Let $c$ be a simple closed curve in $F$.   $2c^\bullet=0$ if $c$ bounds a blank disk with more than one singular point, where by a blank disk we mean a disk that has no curve on it. 
\end{lemma}
{\bf Proof.} easy.
\begin{lemma}
\label{sepagenus}
$2 c^\bullet=0$ if $c$ is a separating simple closed curve in $F$ cutting a genus $k$ empty surface $S$ without a boundary other than $c$ itself. 
\end{lemma}
\noindent
{\bf Proof.} First we may assume that $S$ does not contain any singular point: if it does, get rid of it by the following relation:
$$2\ \raisebox{-.15in}{\psfig{figure=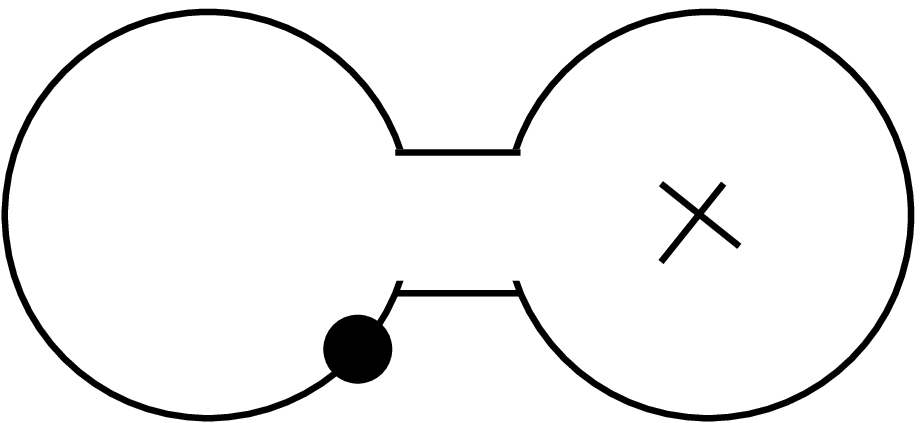,height=1cm}}=\raisebox{-.15in}{\psfig{figure=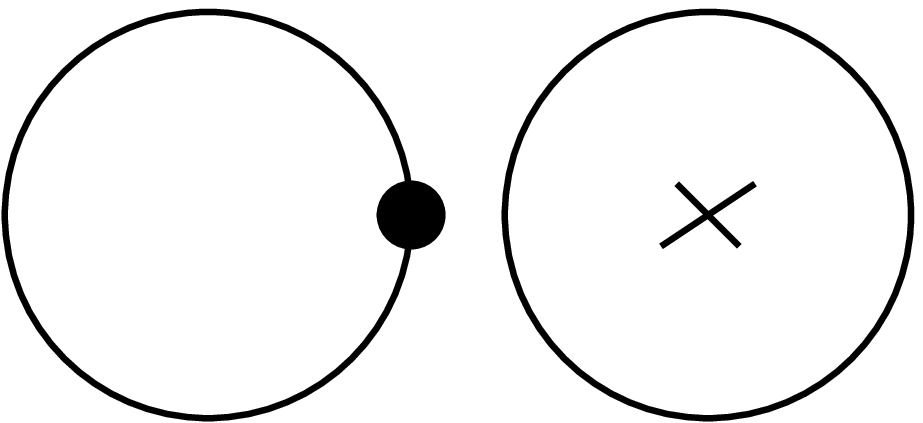,height=1cm}} + \raisebox{-.15in}{\psfig{figure=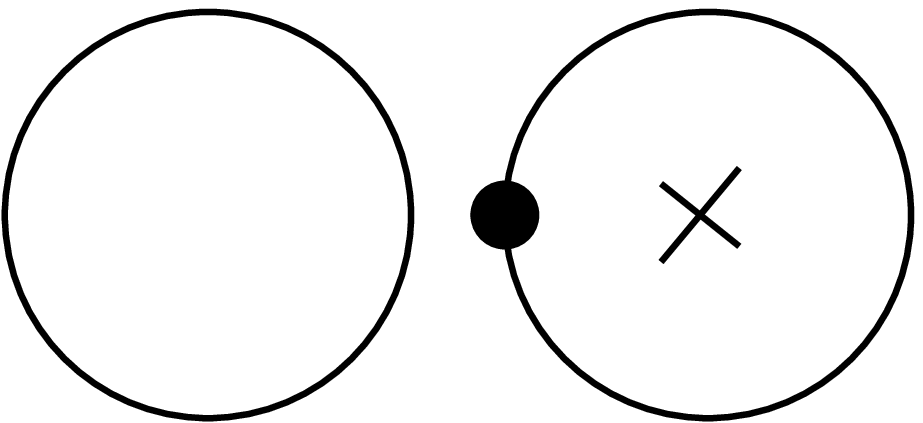,height=1cm}}=2\ \raisebox{-.15in}{ \psfig{figure=survived.eps,height=1cm}} 
 . $$
Apply the annular skein relation to the following series of pictures appropriately. \\
$$\raisebox{-.5in}{\psfig{figure=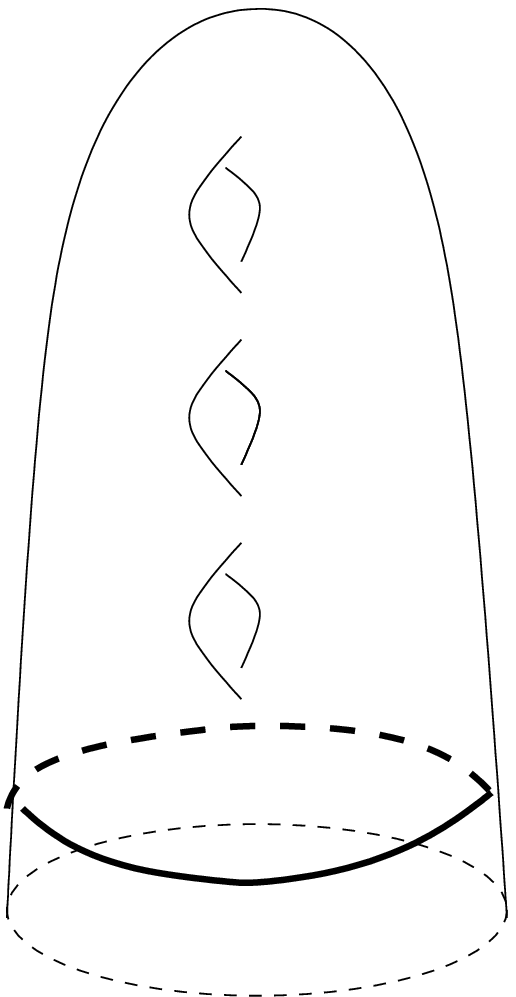,height=3cm}}  \to \raisebox{-.5in}{ \psfig{figure=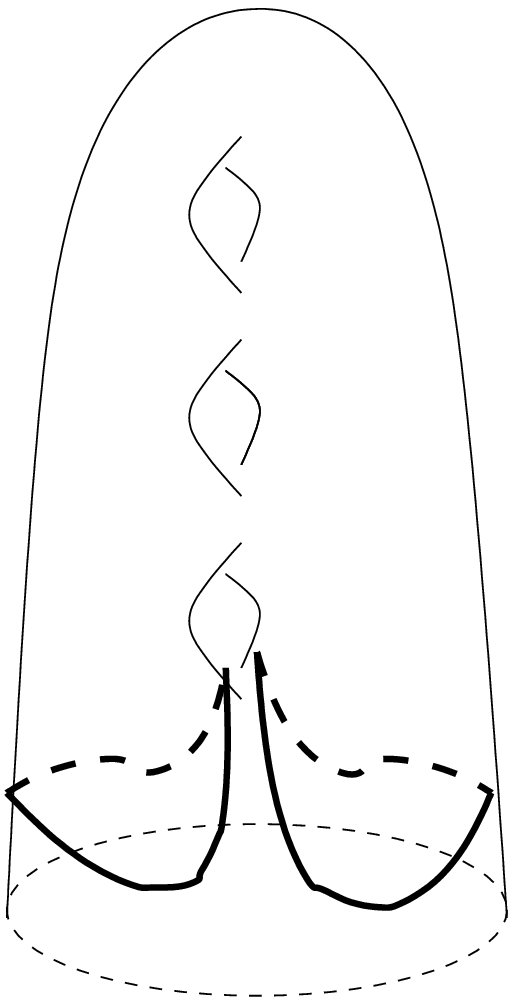,height=3cm}} 
\to \raisebox{-.5in}{\psfig{figure=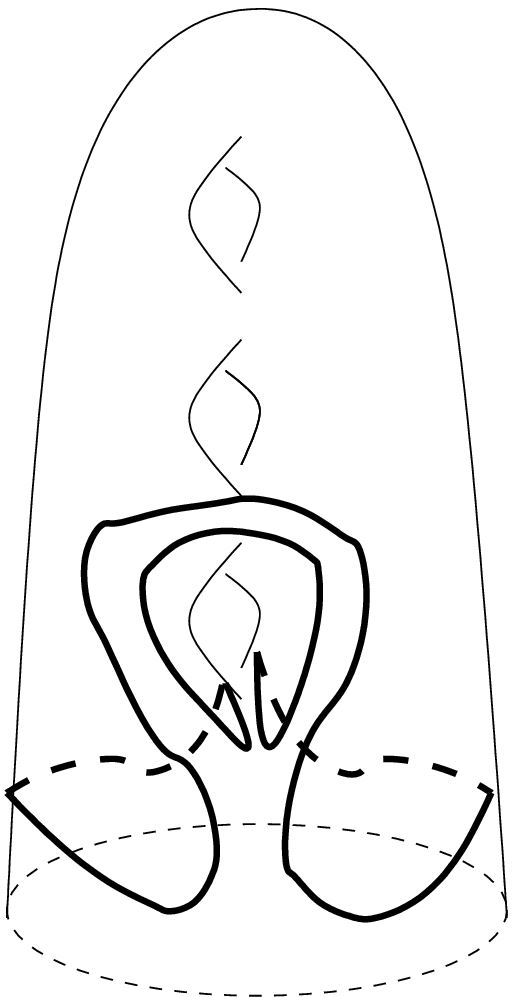,height=3cm}} \to \raisebox{-.5in}{\psfig{figure=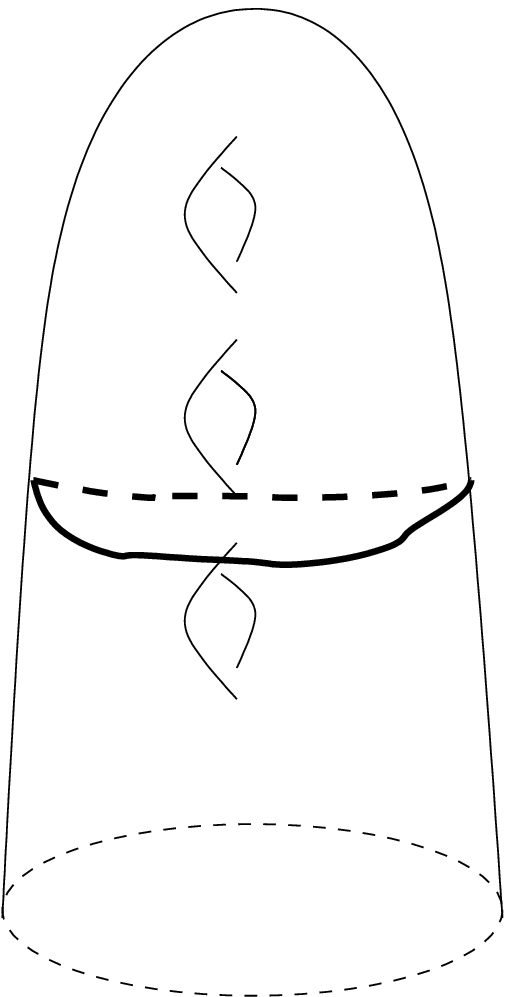,height=3cm}}.  $$
Repeating the process $k$ times one obtains
$ 2c^\bullet = 2c'^\bullet, $ where $c'$ bounds a disk, and thus the right hand side is zero.  \qed.  \\ \ \\
%%%%%%% arc things exported %%%%%
\begin{lemma}
\label{sliding}
Let $\al$ and $\al'$ be states with a dotted circle $c^\bullet$, $c'^\bullet$
respectively. Assume all the other components are the same, and that $c$ and $c'$ bound a surface which is blank in both states. Then $2\al=2\al'$. 
\end{lemma}
{\bf Proof.} similar to the previous lemma. 
\begin{lemma}
\label{switch1}
Let $ \al, \al' \in SBN(F,R)$ be states consisting of the same curves, with a dot in each of two circles cobounding an blank surface without other boundaries. Then $\al =-\al'$. 
\end{lemma}
\noindent
{\bf Proof.} \\
Using the previous lemma, we have
$$\raisebox{-.15in}{\psfig{figure=2curves1.eps,height=1cm}} + \raisebox{-.15in}{\psfig{figure=2curves2.eps,height=1cm}}=2\ \raisebox{-.15in}{\psfig{figure=glued.eps,height=1cm}}=0,$$
where the "clouds" in the circles signify that there may be something there, such as boundaries, singularities, genus, rather than blank disk.  \qed
%%%%%%%%%%%%
\begin{lemma} 
A state $\al$ with more than one dot is zero if there is a pair of dotted circles that are connected by a path of circles. 
\end{lemma}
{\bf Proof.} 
Let $c$, $c'$ be circles in $\al$ that have a dot on each, and let $\{c_1...c_n\}$ be a "path" of circles connecting $c$ and $c'$, namely  $c$ and $c_1$ are adjacent, $c_i$ and $c_{i+1}$ are adjacent, $c_n$ and $c'$ are adjacent. WLOG one may assume that $\{c_i\}$ and $c$, $c'$ are all distinct, and that $\{c_i\}$ are all white. Let us call $n=d(c,c')$ the distance between the dotted circles $c$ and $c'$ along the path. We show by induction in $n$ that $\al=0$. 
When $n=1$, we have
\begin{eqnarray*}
&&c^{\bullet} \cup c_1  \cup c'^\bullet = -c \cup c^{\bullet}_1  \cup c'^\bullet +2 (c \sharp c_1)^{\bullet}  \cup c'^\bullet = \\
&& -\{-c \cup c_1  \cup c'^{\bullet \bullet} + 2 c \cup (c_1 \sharp c')^{\bullet \bullet} \} 
+2 \{ -(c \sharp c_1) \cup c'^{\bullet \bullet} +2  (c \sharp c_1 \sharp c')^{\bullet \bullet} \}=0. 
\end{eqnarray*} 
For general $n$, we have
\begin{eqnarray*}
&& c^{\bullet} \cup c_1 \cup c_2 \cup \cdots \cup c_n \cup c'^\bullet  = \\
 && -c \cup c^{\bullet}_1 \cup c_2 \cup \cdots \cup c_n \cup c'^\bullet + 2 (c \sharp c_1)^{\bullet} \cup c_2 \cup \cdots \cup c_n \cup c'^\bullet=0,
 \end{eqnarray*}
 noting that $d(c_1, c')=d((c \sharp c_1), c')=n-1$ and applying induction.  \qed. \\
 
 For a state $\al \in SBN(M,R)$, we denote by $|\al|$ the number of circles contained in $\al$. We introduce the filtration and grading for $\al \in SBN(M,R)$ as in section \ref{def}. Let \\
$S_i:= \spa_R \{\mbox{ {\rm state }} \in SBN(F,R) \;  | \;  |\al | \leq i, \mbox{\rm a component is dotted} \}$
  Then we have $$S_0 \subset S_1 \subset S_2 \subset \cdots.$$
%and each inclusion is a submodule with embedding given by $$\al \in S_n \to  \frac{1}{2} \al \sqcup \psfig{figure=circle.eps,height=0.5cm} \in S_{n+1}. $$ 
%$$\psfig{figure=dottedarc.eps,height=1.5cm} =\frac{1}{2} ( \psfig{figure=connected1.eps,height=1.5cm} +\psfig{figure=connected2.eps,height=1.5cm}) ,$$ 
Clearly $\mathop{\lim}_{\to}=SBN(F,R).$\\ Let $G_i=S_i/S_{i-1}$. Since $R$ is a direct product of simple rings, a basis of $G_i$'s give a basis for $SBN(F,R)$. 
 Note that in $G_i$ we have the following relation between dotted circles. 
$$ \label{jump}   \raisebox{-.15in}{\psfig{figure=2curves1.eps,height=1cm}} =
\raisebox{-.15in}{\psfig{figure=glued.eps,height=1cm}} -
\raisebox{-.15in}{\psfig{figure=2curves2.eps,height=1cm}} \sim -
\raisebox{-.15in}{\psfig{figure=2curves2.eps,height=1cm}}.$$ Thus it is easy
to see that a state with more than one dot is zero: one may bring one dot to
another dotted circle with sign change.  Thus from now on we only consider
"dotted states" as those that have a single dot. For a dotted state, any choice of the dotted circle would give rise to the same state up to sign. \\
For a better description of states, we define the following notion:
\begin{definition}
For a dotted state $\al$, we define a stack $e$  to be a non-empty collection of circles $\{c_i\}$  in $\al$ defined inductively: a circle $c$ is in $e$ if there is a circle $c' \in e$ that cobounds a blank surface with $c$.  Clearly every circle in $\al$ belong to some stack. All the circles in a stack are homologous, thus a stack determines an element of $H_1(F,\Z_2)$, which we also denote by $e$, abusing the notation. By $w_e$ we denote the number of the circles $|e|$ in a stack $e$, and call it weight. Note that two distinct stacks may determine the same element of $H_1(F,\Z_2)$.  \\
We define a "band" $N(e)$ for a stack $e$ as follows: let us order the circles in $e$ as $c_0, c_2, ... c_n$, where $n=|e|$, $c_i$ and $c_{i+1}$ are adjacent for $i=0,...,n$. We define $N(e)$ to be a connected surface in $F$ so that $\partial N(e)=c_0 \cup c_n$, and $c_i \subset N(e)$ for all $i$. All $N(e)$'s are disjoint, regardless of choice of $N(e)$ for each $e$: $N(e)$ depends on the ordering of the circles in $e$. 
\end{definition}
\psfig{figure=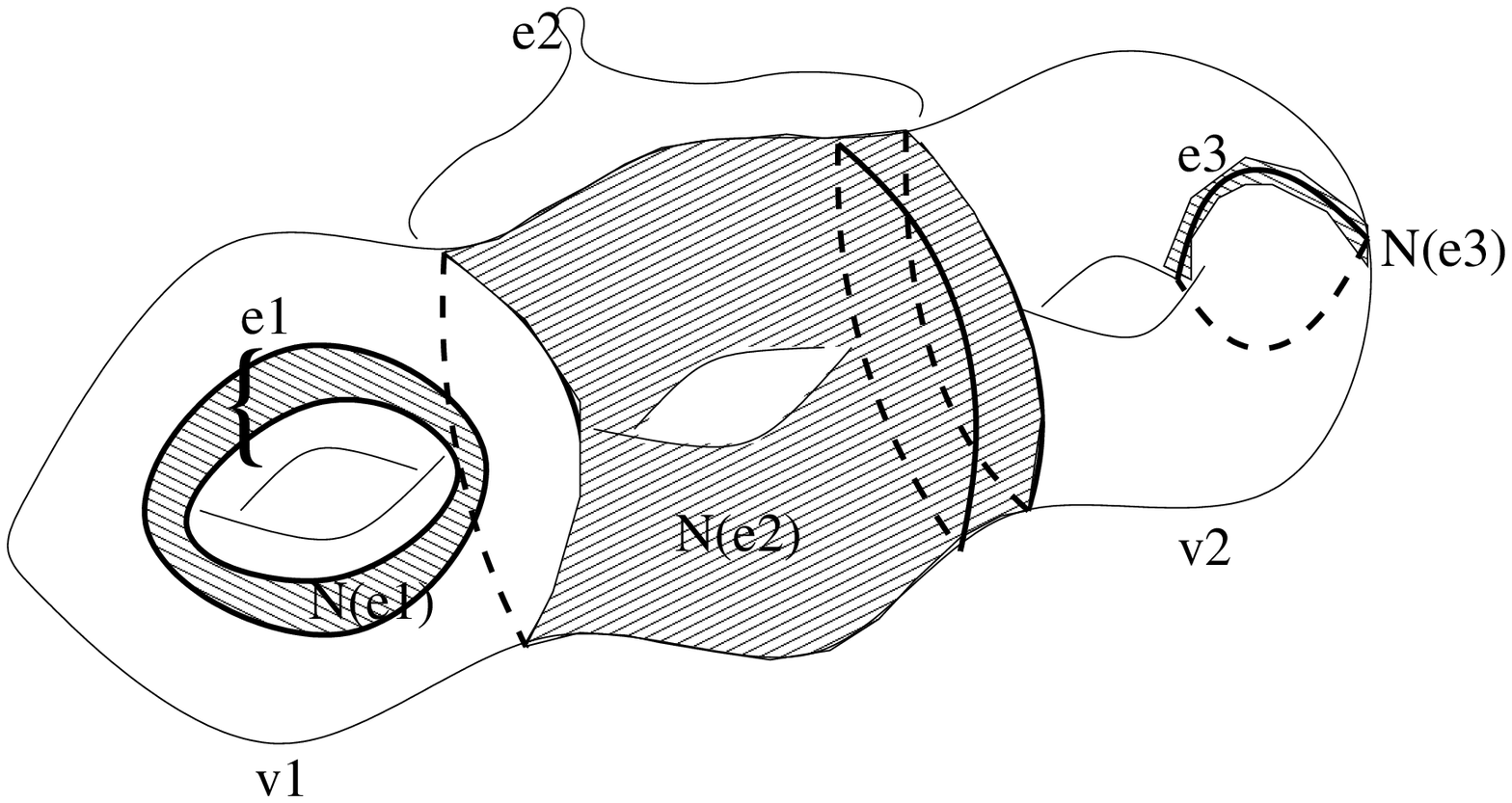,height=2.5cm}  \hspace{1cm}
\psfig{figure=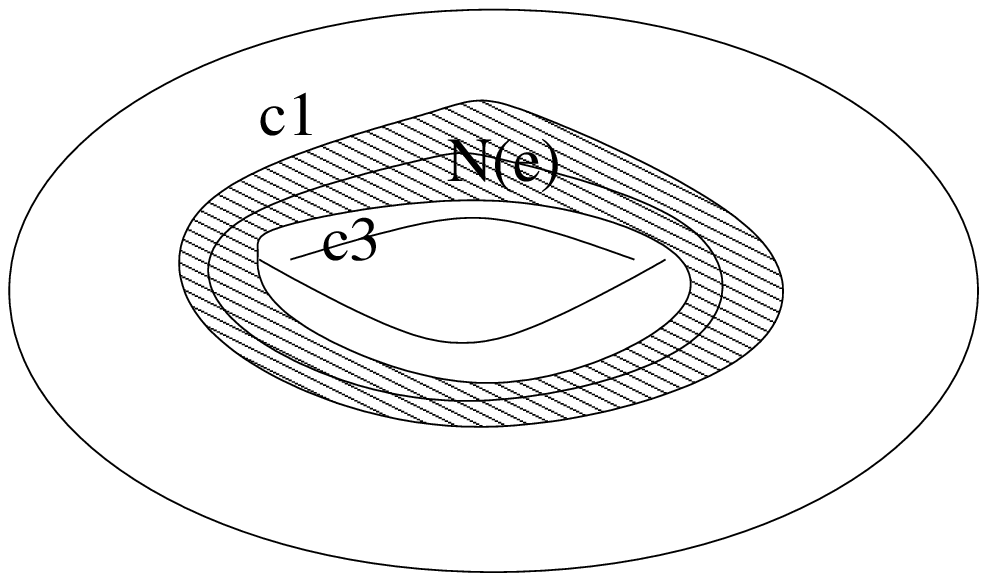,height=2cm} 
\ \\ \ \\
It is easy to see, by using Lemma \ref{sliding} and the fact that one may shift the dot by sign change, that one may move the circles within each stack without changing the state, as long as they do not interfere with circles in other stack. \\
Further, we have the following lemmas:
\begin{lemma}
\label{trivalent}
A dotted state $\al$ is zero  if there is a connected component $F' \subset  F\backslash \sqcup_e \stackrel{\circ}{N(e})$ that contains more than two distinct circles in $\al$ among its boundary components. 
\end{lemma}
{\bf Proof.} \\
Let $c_1$, $c_2$, $c_3$ be distinct circles in $\al$ that are among the boundary component of $F'$, and by $\al_i$ I denote the state with a dot on $c_i$ and has identical circles as $\al$. Then, using the formula \ref{jump} we have
$\al_1=-\al_2=\al_3=-\al_1$, thus $\al=0.$ \qed
\begin{lemma}
\label{loop}
A state $\al$ is zero if there is a connected component $F' \subset  F\backslash \sqcup_e \stackrel{\circ}{N(e})$ and a stack $e$ with $w_e >1$ odd, so that $\partial N(e) \subset F'.$ 
\end{lemma}
{\bf Proof.} This is easily checked, considering the sign changes caused by the moves of the dot. \qed
\\
These lemmas give a strong restriction on the possible non-zero dotted states. In order to enumerate them better, we introduce the following diagram along with some data for each state. 
\begin{definition}(diagram)
For a  state $\alpha \in SBN(F,R)$, the diagram $\Gamma$ of $\al$ consists of the following: \\
\begin{itemize}
\item $\Gamma^{(1)}=\{{\rm edges}\}$. Each  edge $e$  corresponds to a stack $e$. 
\item $\Gamma^{(0)}=\{{\rm vertices}\}$. Each vertex $v$ corresponds to a connected component $F_v$ of $F^o$, where $F^o=F \backslash \sqcup_e \stackrel{\circ}{N(e})$. 
\end{itemize}
An edge $e$ is connected to a vertex $v$ by one end if $\partial F_v \cap N(e)$ is a circle in $e$. $e$ forms a loop around $v$ if $\partial F_v \cap N(e)$ is two circles in $e$. 
%For a graph we have the following data: \\
For $e \in \Gamma^{(1)}$, we have
\begin{itemize}
\item   the weight $w_e=|e|$ is the number of circles that belong to $e$. 
\item the $Z_2$ homology class that each element in $e$ determines. We denote it by $e$, abusing the notation. Note however that it is possible that distinct edges determine the same $Z_2$-homology class. 
\end{itemize}
\end{definition}
An example of a pure state and the corresponding graph is illustrated below:
a pure state $\al$ (the dot is omitted): \\
\psfig{figure=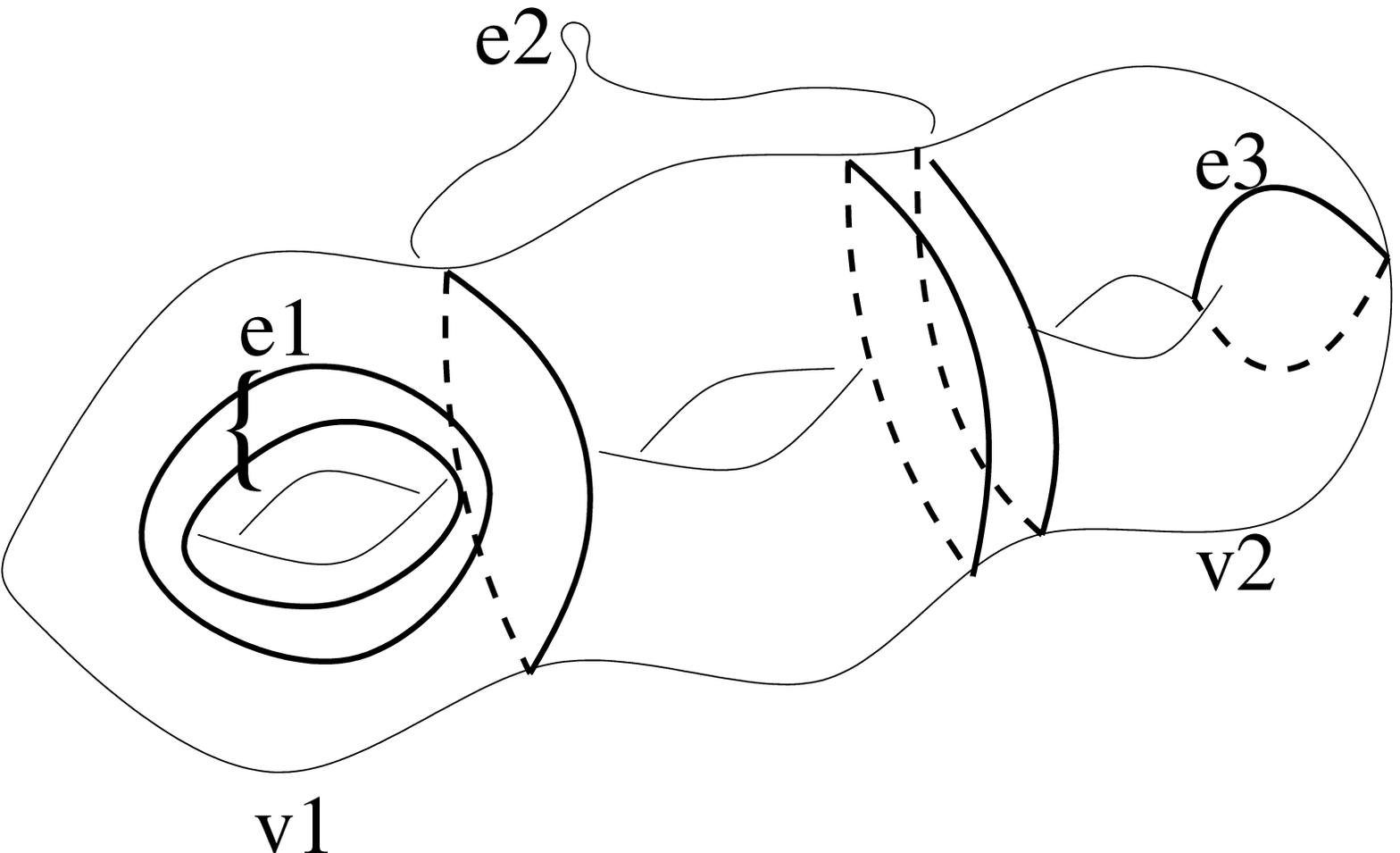,height=2cm} \\
the graph $\Gamma$ of $\al$: \\
\psfig{figure=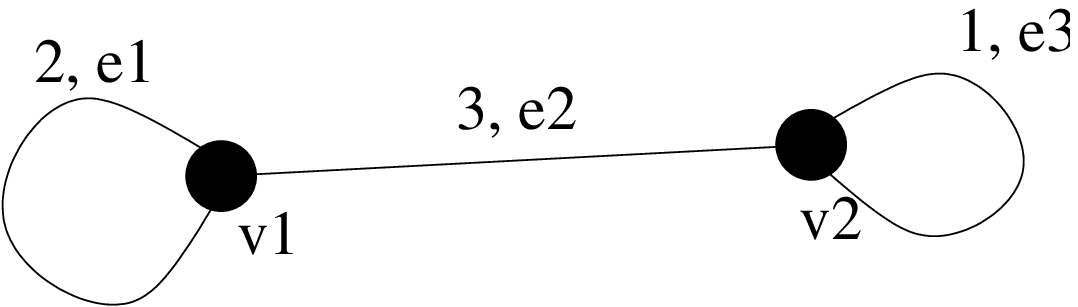,height=1.5cm}\\
\begin{lemma}
For a given surface $F$, the graph along with the data for each edge determines a pure state in $SBN(F,R)$ if for each $v \in \Gamma^{(0)}$ $\sum_{v \in e\neq {\rm loop}} e=0 \in H_1(F,\Z_2)$ holds. 
\end{lemma}
We call a graph that corresponds to a pure state {\it admissible}. \\
{\bf Proof.}\\
First, a graph corresponding to a pure state is admissible since  a vertex $v \in \Gamma^{(0)}$ corresponds to a surface with boundary components specified by the edges that connect to $v$, where a loop will give two circles, thus the contribution is zero in $\Z_2$-homology. On the other hand, for a given admissible graph with the data, the circles corresponding to the edges attached to a vertex must bound a surface in $F$ by the admissibility condition. One many randomly choose circles for each edge as many as the weights. \qed.
\begin{lemma}
For an admissible graph, we have the following reduction rule:
  $$  \raisebox{-.20in}{\psfig{figure=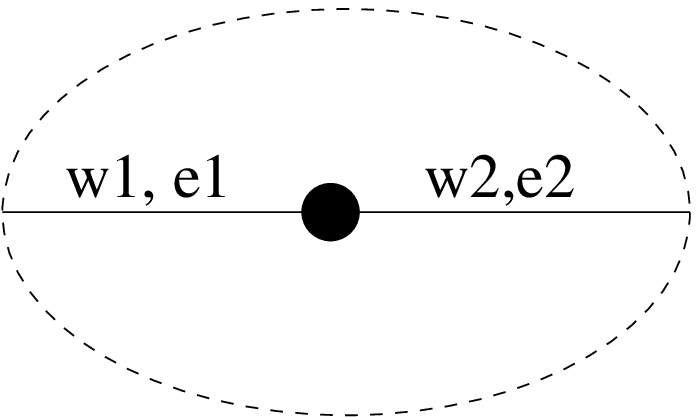,height=1.5cm}}=  \raisebox{-.20in}{\psfig{figure=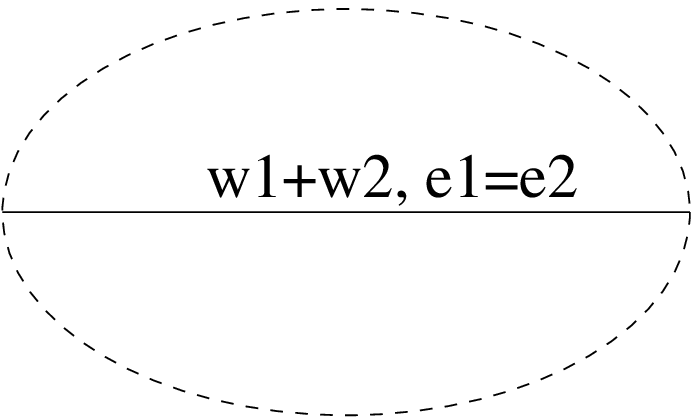,height=1.5cm}}$$
%\item[(2)]  $\psfig{figure=stick.eps,height=1cm}=0$
 
\end{lemma}
{\bf Proof.}
The relation (1) is given by the fact that, by admissibility $e_1=e_2$ in $\Z_2$ homology, and thus the circles in $e_1$ and $e_2$ belong to the same edge. \qed

We call an admissible graph {\it reduced} if there is no divalent vertex connecting two edges. 

An admissible graph $\Gamma$ determines a pure state in the filtration $S_n$ if $|\Gamma|:=\sum_{e \in  \Gamma^{(1)}} w_e \leq n$. 
For the sake of notation economy, we consider admissible graphs  as if they are pure states in $SBN(F,R)$. The Lemmas \ref{trivalent}, \ref{loop} lead to the following formulas for the graphs. 
%%%%%%%%
\begin{lemma}
A reduced admissible graph $\Gamma \in $ determines zero in  $G_n=S_n/S_{n-1}$ if:
\begin{enumerate}
\vspace{.25in}
\item  $|\Gamma| \leq n-1$
\vspace{.25in}
\item  $ \raisebox{-.15in}{\psfig{figure=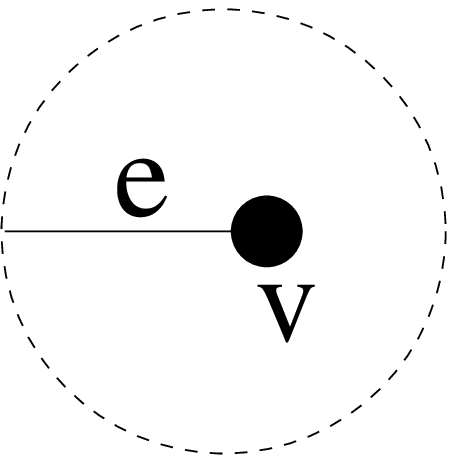,height=1cm}} \subset \Gamma$
\vspace{.25in}
\item for some $v \in \Gamma^{(0)}$,  the number of edges connected to $v$ is more than two. 
\vspace{.25in}
\item  $ \raisebox{-.15in}{\psfig{figure=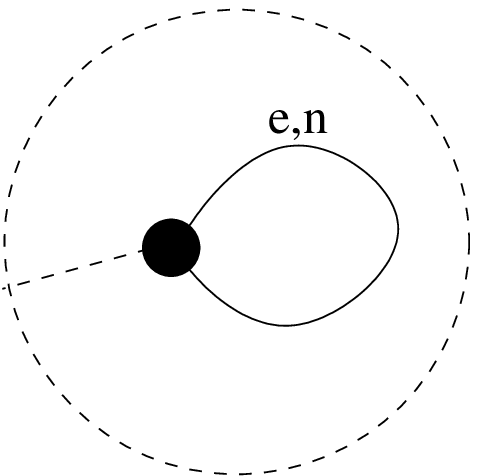,height=1cm}} \subset \Gamma,$ $n>1$, $n$ is odd. 
\vspace{.25in}
%\item   There is a cycle $e_1 e_2 .... e_n$ for $n \geq 2$ that $ \sum_i w_{ei}$ is odd. 
%\vspace{.25in}
%\marginpar{(4) may not be needed.}
%\item[(5)] $\sum_{v \in e} w_e$ is odd for some $v \in \Gamma^{(0)}$
\item  $ \raisebox{-.15in}{\psfig{figure=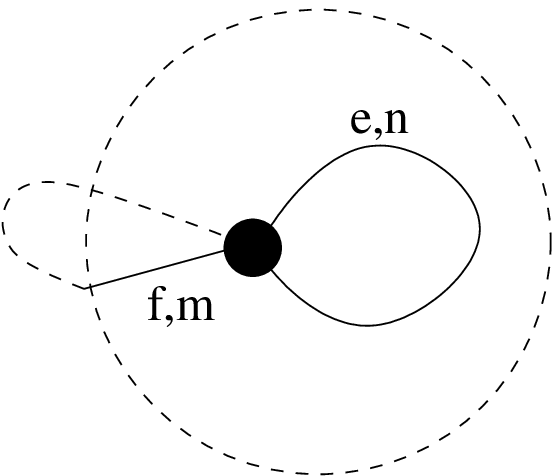,height=1cm}} \subset \Gamma,$ $n$ even, $m \neq 0$
\vspace{.25in}
\end{enumerate}
where dotted line may or may not exist. 
\end{lemma}
%{\bf Proof.} \\
% (1): The circle in the edge adjacent to the surface corresponding to $v$ bounds the surface, thus bringing the dot on the circle with adequate sign change we obtain zero. \\
% (2): Assume that there are three edges connected to $v$. Choose circles $c_1$, $c_2$, $c_3$ from each wire. Denote a state with a dot on $c_i$ by $\al_i$. Then $\al_1=-\al_2=\al_3=-\al_1$, thus $\al_i=0$. \\
% (3) Let $c_1...c_n$ be all the circles in the wire that belong to $e$, where $n$ is odd and $c_{i}$ and $c_{i+1}$ are adjacent. We denote a state with a dot on $c_i$ by $\al_i$. Then $\al_1=-\al_2=... =-\al_{n-1}=\al_n=-\al_1$, thus $\al_i=0$.  \\
% (4) Similar to (3) \\
% (5) Same deal.  QED
Using the previous lemmata we obtain the following classification. 
\begin{theorem}
A reduced admissible graph $\Gamma$ in $G_n$ for some $n$  is  one of the following graphs.
%, where homology classes are omitted. 
\begin{itemize}
\vspace{.25in}
\item[(a)]  \raisebox{-.15in}{\psfig{figure=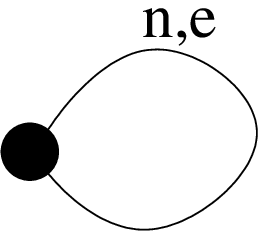,height=1cm}} , $n=1$ or even, $e \in H_1(F,\Z_2)$
\vspace{.25in}
\item[(b)]  \raisebox{-.15in}{\psfig{figure=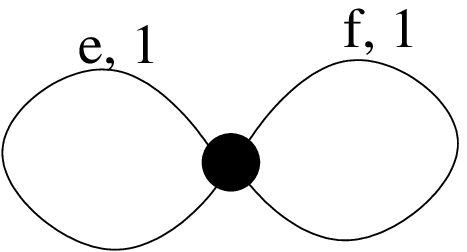,height=1cm}},  $e \neq f \in H_1(F,\Z_2)$
\vspace{.25in}
\item[(c)]   \raisebox{-.15in}{\psfig{figure=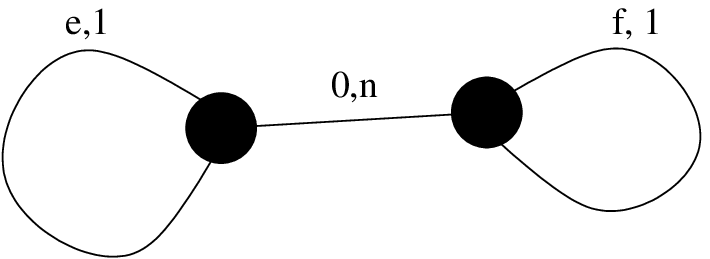,height=1cm}}, $e \neq f \in H_1(F,\Z_2)$. 
\vspace{.25in}
\end{itemize}
\end{theorem}
Note that the graph (b) is the case $n=0$ of the graph (c). 
\begin{cor} 
Let $H_1:=H_1(F,\Z_2)$ and $\Delta=\{ (e,e) \vert e \in H_1\}.$\\
\vspace{.1in} \\
\noindent
 $G_1=RH_1 ={\rm span}_{R} \{ \raisebox{-.15in}{\psfig{figure=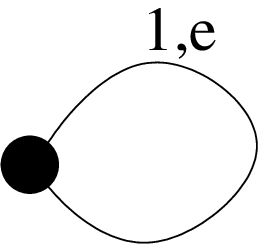,height=1cm}}\vert  e \in H_1 \},  $ \\
\vspace{.1in} \\
%\begin{eqnarray*}
%G_2&=&RH_1(F,\Z_2) \oplus RH_1(F,\Z_2) \times RH_1(F,\Z_2)  \\
%   &=&  {\rm span}_{R} \{\psfig{figure=oneloop2.eps,height=1cm}\vert  e \in H_1(F,\Z_2)\} \oplus 
%   {\rm span}_{R} \{\psfig{figure=oneoneloops.eps,height=1cm}\vert  e,f \in H_1(F,\Z_2)\} 
%   \end{eqnarray*}
 $G_n=$ \\ \ \\
$ \left\{ \begin{array}{cc}  {\rm span}_{R} \{  \raisebox{-.15in}{\psfig{figure=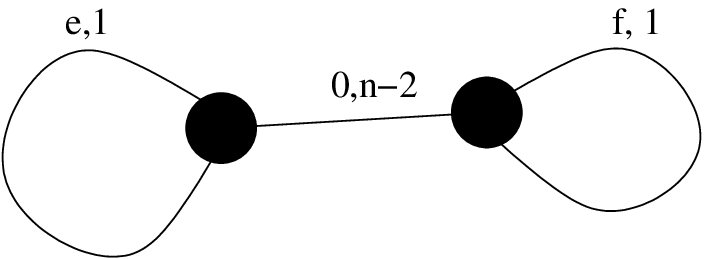,height=1cm}} \vert e \neq f \in H_1   \} \hspace{7cm}   & \\ &\ \\= 
  R(H_1 \times H_1 \backslash \Delta) \quad \mbox{if } n \mbox{ is odd}>1,\hspace{7cm} & \\ 
  & \ \\
  {\rm span}_{R} \{  \raisebox{-.15in}{\psfig{figure=oneloop.eps,height=1cm}} \vert e \in H_1 \} \oplus  \ {\rm span}_{R} \{  \raisebox{-.15in}{\psfig{figure=glasses1.eps,height=1cm}} \vert e \neq f \in H_1   \} \hspace{2.5cm} &
 \\ \ & \\
 = RH_1 \oplus R(H_1 \times H_1 \backslash \Delta)  \quad 
  \mbox{if } n \mbox{ is even}>0 \hspace{7cm} &  \end{array} \right.
  $
 \end{cor}
 \vspace{0.25in}
%%%%%%%%%%
\subsection{The structure of $BN_v(M,R)$ and $BN_h(M,R)$} \ \\ \ \\
 In the previous subsection we clarified the structure of $SBN(F,R)$, which was known to contain $BN_v(M,R)$. In fact we have the following:
  \begin{proposition}
 \label{bnsbn}
 $$BN_v(M,R)=SBN(F,R).$$
 \end{proposition}
 \noindent
 {\bf Proof} \\
We maintain the notation used in the previous subsection regarding $SBN(F,R)$ . Let $F'_n:={\rm span}_R\{ \mbox{ states } \alpha \in BN_v(M,R) | \ |\al| \leq i \}$, and $G'_n=F'_n \slash F'_{n-1}$. Since we already know that states consisting of incompressible surfaces without dots are linearly independent, we show that  the preimages of the basis $ \{\al_i\}$  of $G_n$ with respect to the fibration map $\psi$ are linearly independent in $G'_n$. Assume that the component that the dot in each $\al_i$ is placed is arbitrarily  chosen and fixed for all. We show that the states $\{\al'_i\}$,  where $\al'_i=\psi^{-1}(\al_i) \cong \al_i \times S^1$ are linearly independent in $G'_n$.   We proceed in a similar way to Proposition \ref{independence}. Let $\varphi_i: BN_v(M,R) \to  \mathbb{R}$ be defined as follows: For a state $\beta$, if it is a disjoint union of $\al'_i$ possibly with the dot dislocated, and singly dotted spheres bounding balls, and $k$ compressible white tori,  let $\varphi_i(\beta)=(-1)^l (1/2)^k$, where $l$ is the number of shifts that the dot in $\beta$ makes to arrive to the default position of the dot in $\alpha'_i$ within the surfaces corresponding to $\al'_i$. Note that shifting of a dot among the components of the surface $\al'_i$ is fully expressed in terms of components of the curve  in $\al_i$, which has been discussed earlier.  Since $\al_i$ is not zero, $l$ is well-defined modulo $2$.  Now we check that $\varphi_i$ is compatible with the BN relations modulo grading. The only non-trivial relation to check is (2). Let 
\begin{center} 
$\al'_\infty=$\raisebox{-.25in}{$ \psfig{figure=rel2left.eps,height=1.5cm}$}, \
 $\al'_+=$\raisebox{-.25in}{$ \psfig{figure=rel2right1.eps,height=1.5cm}$} , \
$\al'_-=$\raisebox{-.25in}{$\psfig{figure=rel2right2.eps,height=1.5cm} $},
\end{center}
where by these equalities we denote pure states with local modifications illustrated by the right hand sides. We need to show that $\varphi_i(\al'_\infty)=\varphi_i(\al'_+) + \varphi_i(\al'_-)$. We may assume that $\al'_\pm$ does not have any spheres bounding balls, or tori bounding solid tori.   If the underlining surfaces of $\al'_\pm$ is compressible, or incompressible but not isotopic to that of $\al'_i$, both side of the equality is zero. Otherwise $\al'_+=\pm \al_i=-\al'_i$, thus we have
$$ \varphi_i(\al'_+) + \varphi_i(\al'_-)=(-1)^{l} \varphi_i(\al'_i)+(-1)^{l+1} \varphi_i(\al'_i)=0.$$
 The number of the components of $\al'_\infty$ is $n-1$, thus $\varphi_i(\al'_\infty)=0$ in $G'_n$. Thus we proved the independence of  $\{\al'_i\}$ in $G'_n$, which implies $BN_v(M,R) \cong SBN(F,R)$. 
 \\ \ \\ \indent
 Next we want to understand the horizontal surfaces.  If $M$ has a horizontal surface, let $H$ be a connected horizontal surface so that $\psi|_H:H \rightarrow F$. Orienting the fibers of $M$ we can define a map $\phi:H \rightarrow H$ by following the fibers in the positive direction. The map $\phi$ has finite order and $H\times [0,1]/\sim$, where we identify $(x,1)$ with $(\phi(x),0),$ is homeomorphic to $M$. Let
 $H_1(H)^{\phi}$ denote the invariant part of $H_1(H)$ under the action of $\phi$.
 For the sake of computation we orient $H$ and $M$ and $F$.  If $H'$ is any horizontal surface the degree of $\psi_{H'}$ is an integer $d(H')$. Making $H'$ transverse to  $H$ we can take their oriented intersection and view it as defining an element $z(H')$ of $H_1(H)^{\phi}$. Hence we can assign to any connected oriented horizontal surface the ordered pair $(z(H'),d(H')) \in H_1(H)^{\phi}\oplus \mathbb{Z}$. The image of this map is all indivisible elements, which we denote $I(H)$. Two connected oriented surfaces are isotopic if and only if they have the same image in $I(H)$. Let $P(I(H))$ be all elements of $I(H)$ with 
 $d(H')>0$, these points are in one to one correspondence with the isotopy classes of unoriented horizontal surfaces in $M$. If $f \in P(I(H))$, let $f^k$ denote $k$ parallel copies of $f$ without dots. Let $d_f$ denote on copy of $f$ with a dot.
 \begin{theorem}
 
 $BN_h(M,R)= \oplus_{f\in P(I(H))} \{f^k\}_k \oplus \{d_f\}.$ 
 \end{theorem}
 
 The proof is just like the example $S^1\times S^2$, noting that two non-parallel horizontal surfaces cannot be disjoint. 
 
%%%%%%%%%%%%%%%%
% \section{misc}
% mention graph manifolds. 

The authors would like to thank Mikhail Khovanov and Walter Neumann for helpful suggestions. 

\end{document}